\documentclass[a4paper,reqno]{amsart}
\usepackage{hyperref}
\usepackage{enumerate}

\usepackage{a4wide,amsthm}

\usepackage{amsmath}
\usepackage{comment}
\usepackage{amsfonts}
\usepackage{color}

\usepackage{graphicx}
\usepackage{subfig}
\usepackage{bm}
\usepackage{tikz}
\usepackage{amssymb}

\newtheorem{lemma}{Lemma}[section]

\newtheorem{theorem}[lemma]{Theorem}

\theoremstyle{definition}

\theoremstyle{remark}
\newtheorem{remark}[lemma]{Remark}

\newcommand{\dst}{\displaystyle}
\newcommand{\ffi}{\varphi}

\def\A{{\mathbb{A}}}

\def\C{{\mathbb{C}}}

\def\N{{\mathbb{N}}}

\def\R{{\mathbb{R}}}
\def\T{{\mathbb{T}}}

\def\Z{{\mathbb{Z}}}

\def\nn{{\mathcal N}}

\newcommand{\normL}[1]{{\left\|{#1}\right\|_{L^1_K}}}

\newcommand{\norm}[1]{{\left\|{#1}\right\|}}

\newcommand{\abs}[1]{{\left|{#1}\right|}}

\newcommand{\ud}{\,\mathrm{d}}

\usepackage{bbm}
\usepackage{mathtools}
\numberwithin{equation}{section}
\usepackage{graphicx} 

\date{08/03/2024}

\begin{document}
\title[On $L^1$-norms for non-harmonic trigonometric polynomials]{On $L^1$-norms for non-harmonic trigonometric polynomials with sparse frequencies}
\author{Philippe Jaming, Karim Kellay, Chadi Saba \& Yunlei Wang}
\address{Univ. Bordeaux, CNRS, Bordeaux INP, IMB, UMR 5251,  F-33400, Talence, France}



\dedicatory{In honnor of Karlheinz Gr\"ochenig on the occasion of his 60$^{\mbox{th}}$ birthday.}

\begin{abstract}
In this paper we show that, if an increasing sequence $\Lambda=(\lambda_k)_{k\in\Z}$
has gaps going to infinity $\lambda_{k+1}-\lambda_k\to +\infty$
when $k\to\pm\infty$, then for every $T>0$
and every sequence $(a_k)_{k\in\Z}$ and every $N\geq 1$,
$$
A\sum_{k=0}^N\frac{|a_k|}{1+k}\leq\frac{1}{T}\int_{-T/2}^{T/2}
\abs{\sum_{k=0}^Na_ke^{2i\pi\lambda_k t}}\,\mbox{d}t
$$
further, if $\dst\sum_{k\in\Z}\dfrac{1}{1+|\lambda_k|}<+\infty$,
$$
B\max_{|k|\leq N}|a_k|\leq\frac{1}{T}\int_{-T/2}^{T/2}
\abs{\sum_{k=-N}^Na_ke^{2i\pi\lambda_k t}}\,\mbox{d}t
$$
where $A,B$ are constants that depend on $T$ and $\Lambda$ only.

The first inequality was obtained by Nazarov for $T>1$ and the 
second one by Ingham for $T\geq 1$ under the condition that 
$\lambda_{k+1}-\lambda_k\geq 1$. The main novelty is that if those
gaps go to infinity, then $T$ can be taken arbitrarily small.
The result is new even when the $\lambda_k$'s are integers where
it extends a result of McGehee, Pigno and Smith.

The results are then applied to observability of
Schr\"odinger equations with moving sensors.
\end{abstract}

\maketitle 

\section{Introduction}

The aim of this paper is to establish a lower bound of $L^1$-norms of non-harmonic trigonometric polynomials with
sparse frequencies. The results are then applied to obtain $L^1$-observability estimate
of certain PDEs, including the free Schr\"odinger equation. We thus
obtain $L^1$-analogues of a result of Kahane \cite{Ka} and Haraux \cite{Ha} on the $L^2$-norm
of sparse trigonometric polynomials while the $L^2$-observability result was previously obtained
by the first author together with Komornik \cite{JK}.

Let us now be more precise. We first describe the well-known results in the $L^2$-setting.
The celebrated Ingham Inequality gives a lower and upper bound of $L^2([-T,T])$-norms
of (non-harmonic) trigonometric polynomials and is stated as follows:

\begin{theorem}[Ingham \cite{In}]
Let $\gamma>0$ and $T>\dfrac{1}{\gamma}$. Then there
exist constants $0<A_2(T,\gamma)\leq B_2(T,\gamma)$ such that

-- for every sequence of real numbers  $\Lambda=\{\lambda_k\}_{k\in\Z}$ such that $\lambda_{k+1}-\lambda_k\geq\gamma$;

-- for every sequence $(a_k)_{k\in\Z}\in\ell^2(\Z,\C)$,
\begin{equation}
\label{eq:inghamtg}
A_2(T,\gamma)\sum_{k\in\Z}|a_k|^2\leq \frac{1}{T}\int_{-T/2}^{T/2}\abs{\sum_{k\in\Z}a_ke^{2i\pi\lambda_k t}}^2\,\mbox{d}t
\leq B_2(T,\gamma)\sum_{k\in\Z}|a_k|^2
\end{equation}
\end{theorem}

Note that $A_2(T,\gamma),B_2(T,\gamma)$ are explicit constants ({\it see} \cite{KL,JS}).
Ingham has also shown that the upper bound is valid for any $T>0$ while the lower bound may not be true for 
$T\leq \dfrac{1}{\gamma}$. In his seminal work on almost periodic functions \cite{Ka}, Kahane has shown that this condition can be lifted if $\lambda_{k+1}-\lambda_k\to+\infty$ when $k\to\pm\infty$:

\begin{theorem}[Kahane]\label{th:Ka}
Let $\Lambda=\{\lambda_k\}_{k\in\Z}$ such that $\lambda_{k+1}-\lambda_k\to+\infty$
when $k\to\pm\infty$. Then, for every $T>0$,
there exist constants $0<A_2(T,\Lambda)\leq B_2(T,\Lambda)$
such that 
$$
A_2(T,\Lambda)\sum_{k\in\Z}|a_k|^2\leq \frac{1}{T}\int_{-T/2}^{T/2}\abs{\sum_{k\in\Z}a_ke^{2i\pi\lambda_k t}}^2\,\mbox{d}t
\leq B_2(T,\Lambda)\sum_{k\in\Z}|a_k|^2
$$
holds for every sequence $(a_k)_{k\in\Z}\in\ell^2(\Z,\C)$.
\end{theorem}

The constants are not explicit in \cite{Ka}, they were later obtained by Haraux \cite{Ha} (but with constants
that are difficult to compute explicitly, {\it see e.g.} \cite {KL,JS}).

Those inequalities have found many applications in control theory. Among the numerous results
({\it see} the book \cite{KL} for a good introduction to the subject), our starting point
is a result of the first author with V. Komornik \cite{JK}. To state it,
let us introduce some notation. We write $\T=\R/\Z$ and $H^2(\T)=\{f \in L^2(\T) : \sum_{k\in \Z} (1+|k|^2)^2 |c_k(f)|^2 < \infty   \}$, where the $c_k(f)$'s are the Fourier coefficient of $f$. Then the following holds:

\begin{theorem}\label{th:JK}
Fix $(t_0,x_0)\in\R^2$, $a\in\R$ and $T>0$. 
For $u_0\in H^2(\T)$, let $u$ be the solution of
$$
\begin{cases}
u_t=\dfrac{i}{2\pi}u_{xx}&\mbox{in }\R\times\T,\\
u(0,x)=u_0(x)&\mbox{for }x\in \T.
\end{cases}
$$
\begin{enumerate}
\renewcommand{\theenumi}{\roman{enumi}}
\item There exists $D_2(T,a)$ such that, for every $u_0\in H^2(\T)$,
\begin{equation*}
\int_0^T\abs{u(t_0+t,x_0+at)}^2\,\mbox{d}t\leq D_2(T,a)\|u_0\|_{L^2}^2
\end{equation*}
\item If $a\notin\Z$, then there exists $C_2(T,a)$ such that, for every $u_0\in H^2(\T)$,
\begin{equation}\label{33}
C_2(T,a)\|u_0\|_{L^2}^2 \leq \int_0^T\abs{u(t_0+t,x_0+at)}^2\,\mbox{d}t 
\end{equation}
also holds.

\item If $a\in\Z$, then there exists $u_0\not=0$ such that $u(t_0+t,x_0+at)=0$ so that
\eqref{33} fails.
\end{enumerate}
\end{theorem}

Let us sketch the proof. If we write $\dst
u_0(x)=\sum_{k\in \Z}c_ke^{2i\pi k x}
$ then the solution of the Schr\"odinger equation can be written as a Fourier series
$\dst u(t,x)=\sum_{k\in \Z}c_ke^{2i\pi (k^2 t+k x)}$ and the fact that $u_0\in  H^2(\T)$
implies that this series is uniformly convergent. One can thus restrict it to a segment:
$$
u(t_0+t,x_0+2at)=\sum_{k\in \Z}c_ke^{2i\pi k^2(t_0+t)+2i\pi k(x_0+2at)}
:=\sum_{k\in \Z}d_ke^{2i\pi \lambda_k t}.
$$
Then one shows that the $\lambda_k$'s are such that Kahane's Theorem applies (provided $a$ is not an integer).
Our aim is to extend this argument to the $L^1$-setting.

The first task is thus to obtain an $L^1$-version of Ingham's Inequality.
An $L^1-\ell^\infty$ estimate was obtained by Ingham \cite{In} (and is an easy adaptation of the $L^2$ proof)
and a much more evolved $L^1$ to weighted $\ell^1$-inequality was obtained by Nazarov, inspired
by the proof of Littlewood's conjecture by McGehee-Pigno-Smith. The results are the following:

\begin{theorem}
\label{th:AA}
Let $(\lambda_k)_{k\in\Z}$ be an increasing sequence of real numbers such that
there exists $\gamma>0$ for which $\lambda_{k+1}-\lambda_k\geq\gamma$ for every $k$.
Let $(a_k)_{k\in\Z}$ be a sequence of complex numbers.

\begin{itemize}
\item {\em Ingham \cite{In2}\,:} For $T\geq \dfrac{1}{\gamma}$, there exists a constant $A_1(T,\gamma)>0$ such that, for every $N\geq 1$,
$$
A_1(T,\gamma)\max_{k=-N,\ldots,N}|a_k|\leq \frac{1}{T}\int_{-T/2}^{T/2}\abs{\sum_{k=-N}^Na_ke^{2i\pi \lambda_kt}}\,\mathrm{d}t.
$$

\item {\em Nazarov \cite{Na}\,:} For $T> \dfrac{1}{\gamma}$, there exists a constant $\tilde{A_1}(T,\gamma)>0$ such that, for every $N\geq 1$,
$$
\tilde A_1(T,\gamma)\sum_{k=0}^N\frac{|a_k|}{1+k}\leq \frac{1}{T}\int_{-T/2}^{T/2}\abs{\sum_{k=0}^Na_ke^{2i\pi \lambda_kt}}\,\mathrm{d}t.
$$
\end{itemize}
\end{theorem}

Ingham established the first inequality for $T> \dfrac{1}{\gamma}$ in \cite{In}
and improved his result in \cite{In2} showing that it holds when $T=\dfrac{1}{\gamma}$, and that one may take $A_1(T,\gamma)=\dfrac{1}{2}$.
This was further improved by Mordell \cite{Mo}.
There is a major difference between the two inequalities: the right hand side in Ingham's Inequality
is generally much smaller than in Nazarov's Inequality ({\it e.g.} take $|a_k|=1$ for all $k$
then Ingham provides a constant lower bound while Nazarov provides a logarithmic one).
On the other hand, in Nazarov's inequality the sum starts at $0$ and may fail for symmetric sums.
Also its validity for $T= \dfrac{1}{\gamma}$ is an open question.

Further Nazarov did not provide an estimate of the constant $C(\gamma,T)$. However, his proof can be modified to establish quantitative bounds. This was done in \cite{JKS} when $\gamma T$ is large enough and in \cite{JS}
for $\gamma T$ near to $1$.

This result is sufficient to partially extend Theorem \ref{th:JK} to the $L^1$-setting.
The only thing that would be missing is that in Theorem \ref{th:JK}, there is no minimal time
needed thanks to Kahane's extension of Ingham's inequality. However, so far this is unknown
in the $L^1$-case and our first result is precisely to prove this:

\begin{theorem}\label{t11}
Let $\Lambda=(\lambda_k)_{k\in\Z}$ be an increasing sequence with $\lambda_{k+1}-\lambda_k\to+\infty$ when 
$k\to\pm\infty$. 
Then, for every $T>0$, there exists a constant $\tilde A_1(T,\Lambda)>0$ such that,
if $(a_k)_{k\in\N}\subset\C$ is a sequence of complex numbers, and $N\geq 1$, then
\begin{equation}
\label{eq:harrauxl1}
\tilde A_1(T,\Lambda)\sum_{k=0}^\infty\frac{|a_k|}{1+k}\leq\frac{1}{T}\int_{-T/2}^{T/2}\abs{\sum_{k=0}^N a_ke^{2i\pi\lambda_k t}}\,\mathrm{d}t.
\end{equation}
If further $\dst\sum_{k\in\Z}\dfrac{1}{1+|\lambda_k|}$ converges, then
there also exists a constant  $A_1(T,\Lambda)$ such that, 
for every $(a_k)_{k\in\Z}\subset\C$ and every $N\geq 1$,
\begin{equation}
\label{eq:inghamharl1}
A_1(T,\Lambda)\max_{k=-N,\ldots,N}|a_k|\leq\frac{1}{T}\int_{-T/2}^{T/2}\abs{\sum_{k=-N}^N a_ke^{2i\pi\lambda_k t}}\,\mathrm{d}t.
\end{equation}
\end{theorem}

The main difficulty in the proof of this result is that both Kahane's and Haraux's argument can not
be adapted directly. Indeed, both use in a crucial way that in Ingham's Inequality the $L^2$
norm of a trigonometric polynomial is both lower and upper bounded by the $\ell^2$-norm
of its coefficients. In the $L^1$-case, the upper bound is in terms of the $\ell^1$-norm
of the coefficients and does not match the lower bound. Instead, our proof uses a compactness argument so that we don't obtain an estimate of $A_1(\gamma,T),\tilde A_1(\gamma,T)$ in this case.
It would be interesting to obtain such an estimate.

Finally, we apply this result to an observability inequality for the Schr\"odinger equation with moving sensor.
We show the following: take $u_0\in H^2(\T)$ and write its Fourier series $u_0(x)=\sum_{k\in\Z}c_ke^{2i\pi kt}$.
Let $u$ be the solution of
$$
\begin{cases}
u_t=\dfrac{i}{2\pi}u_{xx}&\mbox{in }\R\times\T,\\
u(0,x)=u_0(x)&\mbox{for }x\in \T.
\end{cases}
$$
then, for every $a\notin\Z$ and every $T>0$, there exists a constant
$C(a,T)>0$ such that
$$
\frac{1}{T}\int_0^T |u(t_0+t,x_0+at)|\ud t \geqslant C\sum_{k\in \Z}\frac{|c_k|}{1+|k|}.
$$
Similar results are then obtained for higher order Schr\"odinger equations.

The remaining of the paper is organized as follows. In Section 2 we prove Theorem \ref{t11}.
We then devote section 3 to the free Schr\"odinger equation while the last section
is devoted to higher order Schr\"odinger equations.

\section{Proof of Theorem \ref{t11}}

First note that replacing the sequence $(\lambda_k)_{k\in\Z}$ with a translate $\lambda_k+\lambda$,
leaves \eqref{eq:harrauxl1}-\eqref{eq:inghamharl1} unchanged. So there is no loss of generality in assuming
that $\lambda_0>0>\lambda_{-1}$. We now fix $T>0$.

Define $K$ to be an integer such that, if $|k|\geq K$, $\lambda_{k+1}-\lambda_k\geq\dfrac{2}{T}$.
As a consequence, from Nazarov's inequality, the following holds for 
every sequence $(b_k)_{k\in\N}$ and every $N\geq 0$:
\begin{eqnarray}
\frac{1}{T}\int_{-T/2}^{T/2}\left|\sum_{k=K}^{K+N} b_ke^{2i\pi \lambda_k t}\right|\,\mathrm{d}t
&=&\frac{1}{T}\int_{-T/2}^{T/2}\left|\sum_{k=0}^{N} b_{k+K}e^{2i\pi \lambda_{k+\nn} t}\right|\,\mathrm{d}t\notag\\
&\geq& \tilde A_1\left(T,\dfrac{2}{T}\right)\sum_{k=0}^N\frac{|b_{k+K}|}{k+1}
\geq \tilde A_1\left(T,\dfrac{2}{T}\right)\sum_{k=K}^{K+N}\frac{|b_{k}|}{k+1};
\label{eq:HL2}
\end{eqnarray}
while Ingham's inequality shows that
\begin{equation}
\frac{1}{T}\int_{-T/2}^{T/2}\left|\sum_{K\leq|k|\leq K+N} b_ke^{2i\pi \lambda_k t}\right|\,\mathrm{d}t
\geq A_1\left(T,\dfrac{2}{T}\right)\max_{K\leq|k|\leq K+N}|b_k|.
\label{eq:ingh}
\end{equation}

We first prove \eqref{eq:harrauxl1}.
To do so, we will adopt the following convention. An element $z$ of $\C^{N}$ will be indexed starting at $0$,
$z=(z_0,\ldots,z_{N-1})$. We will identify it with a vector in $\C^M$, $M\geq N$ as well as with a sequence
$(z_k)_{k\geqslant 0}$ by adding $0$'s at the end, {\it i.e.} setting $z_k=0$ for $k\geq N$.
An element of $\C^N$ is thus called a vector or a sequence, which ever is the most convenient.

On $\C^{N}$, we introduce two norms through
$$
\norm{(a_0,\ldots,a_{N-1})}_{\ell^{1,-1}_N}=\sum_{k=0}^{N-1}\frac{|a_k|}{1+k}
$$
and
$$
\norm{(a_0,\ldots,a_{N-1})}_{L^1_N}=
\norm{\sum_{k=0}^{N-1}a_ke^{2i\pi\lambda_k t}}_{L^1([-T/2,T/2])}:=
\frac{1}{T}\int_{-T/2}^{T/2}\abs{\sum_{k=0}^{N-1}a_ke^{2i\pi\lambda_k t}}\,\mbox{d}t.
$$
The first one is clearly a norm while for the second one, it is enough to notice that the 
set $\{t\to e^{2i\pi\lambda t}\}_{\lambda\in\R}$ is linearly independent in $L^1([-T/2,T/2])$. 

As
$\norm{\cdot}_{\ell_N^{1,-1}}$ and $\norm{\cdot}_{L^1_N}$ are both norms
on the finite dimensional space $\C^{N}$, they are equivalent. Thus there
are $\kappa_N\leq 1\leq \Lambda_N$ such that, for every $a\in\C^{N}$,
\begin{equation}
\label{eq:equivnorm}
\kappa_N\norm{a}_{\ell_N^{1,-1}}\leq\norm{a}_{L^1_N}\leq\Lambda_N\norm{a}_{\ell_N^{1,-1}}.
\end{equation}
Nazarov's theorem asserts that one may choose $\kappa_N$ independent of $N$ provided $T$ is large enough.
Our aim is to show that this is possible for every $T$ under our additional condition on $(\lambda_k)_{k\geq 0}$.

Assume towards a contradiction that this is not the case. Then, for every integer $n\geq 1$, there exist
an integer $K_n$ and $a^{(n)}=(a_0^{(n)},\ldots,a_{K_n-1}^{(n)})\in\C^{K_n}$ such that
$\|a^{(n)}\|_{\ell^{1,-1}_{K_n}}=1$ while $\|a^{(n)}\|_{L^1_{K_n}}\leq\dfrac{1}{n}$. The first observation is that 
$K_n\to+\infty$
otherwise, we would contradict \eqref{eq:equivnorm} when $n$ is large enough.
Hence, without loss of generality, we will assume that $K_{n+1}>K_n\geq K$ for every $n$,
where $K$ was defined so that if $|k|\geq K$, $\lambda_{k+1}-\lambda_k\geq\dfrac{2}{T}$.

Next, we split $a^{(n)}$ into two vectors
$$
a^{(n)}_-=(a_0^{(n)},\ldots,a_{K-1}^{(n)},0,\ldots,0)
\quad\mbox{and}\quad
a^{(n)}_+=a^{(n)}-a^{(n)}_-.
$$
With an obvious abuse of notation, we consider that $a^{(n)}_-\in\C^{K}$.
In particular $\|a^{(n)}_-\|_{\ell^{1,-1}_{K}}\leq\|a^{(n)}_-\|_{\ell^{1,-1}_{K_n}}\leq 1$.
Thus, up to taking a subsequence, we may assume that $a^{(n)}_-\to(a_0,\ldots,a_{K-1})$.

Next, define the following functions:
\begin{enumerate}
\item the functions $\ffi^{(n)}$ given by
$$
\ffi^{(n)}(t)=\sum_{k=0}^{K_n}a_k^{(n)}e^{2i\pi\lambda_k t}
$$
so that $\|\ffi^{(n)}\|_{L^1([-T/2,T/2])}\leq\dfrac{1}{n}\to 0$ {\it i.e.} $\ffi^{(n)}\to 0$
in $L^1([-T/2,T/2])$.

\item The functions $\ffi^{(n)}_-,\ffi_-$ given by
$$
\ffi^{(n)}_-(t)=\sum_{k=0}^{K-1}a_k^{(n)}e^{2i\pi\lambda_k t}
\quad\mbox{and}\quad
\ffi_-(t)=\sum_{k=0}^{K-1}a_ke^{2i\pi\lambda_k t}.
$$
This functions are in a finite dimensional subspace of $L^1([-T/2,T/2])$
so that the convergence $a^{(n)}_k\to a_k$ for $k=0,\ldots,K-1$ implies that
$\ffi^{(n)}_-\to \ffi_-$ in $L^1([-T/2,T/2])$.

\item The functions 
$$
\ffi^{(n)}_+=\ffi^{(n)}-\ffi^{(n)}_-=\sum_{k=K}^{K_n}a_k^{(n)}e^{2i\pi\lambda_k t}.
$$
\end{enumerate}

Note that $\ffi^{(n)}_+\to -\ffi_-$ in $L^1([-T/2,T/2])$. On the other hand, for $n\geq m$ we can apply
\eqref{eq:HL2} to $\ffi^{(n)}_+-\ffi^{(m)}_+$ leading to
\begin{eqnarray*}
\frac{1}{T}\int_{-T/2}^{T/2}|\ffi^{(n)}_+(t)-\ffi^{(m)}_+(t)|\,\mbox{d}t
&=&
\frac{1}{T}\int_{-T/2}^{T/2}\abs{\sum_{k=K}^{K_n}\bigl(a_k^{(n)}-a_k^{(m)}\bigr)e^{2i\pi\lambda_k t}}\,\mbox{d}t\\
&\geq& \tilde A_1\left(T,\dfrac{2}{T}\right)\sum_{k=K}^{K_n}\frac{\bigl|a_k^{(n)}-a_k^{(m)}\bigr|}{k+1}.
\end{eqnarray*}
Using also that $a_k^{(n)}\to a_k$ for $k=0,\ldots,K-1$ this shows that
$\bigl(a_k^{(n)}\bigr)_{k\geq 0}$ is a Cauchy sequence in the Banach space
$$
\ell^{1,-1}=\left\{(b_k)_{k\geq 0}\,:\|(b_k)\|_{\ell^{1,-1}}:=\sum_{k=0}^{+\infty}\frac{|b_k|}{k+1}\right\}.
$$
In particular, $\bigl(a_k^{(n)}\bigr)_{k\geq 0}\to (a_k)_{k\geq 0}$ in $\ell^{1,-1}$.
This implies that, for all $k$, $a_k^{(n)}\to a_k$ and that
$$
1=\|a^{(n)}\|_{\ell^{1,-1}_{K_n}} =\|a^{(n)}\|_{\ell^{1,-1}}\to \|a\|_{\ell^{1,-1}}.
$$
We will thus reach a contradiction if we show that $a_k=0$ for all $k$.

To do so, we introduce further functions via
$$
\Phi^{(n)}_\pm(x)=\int_0^x\ffi^{(n)}_\pm(t)\,\mbox{d}t
\quad\mbox{and}\quad
\Phi_-(x)=\int_0^x\ffi_-(t)\,\mbox{d}t
=\frac{1}{2i\pi}\sum_{k=0}^{K-1}\frac{a_k}{\lambda_k}\bigl(e^{2i\pi\lambda_kx}-1\bigr).
$$
Note that as $\ffi^{(n)}_\pm\to\pm\ffi_-$ in $L^1([-T/2,T/2])$,
$\Phi^{(n)}_\pm\to\pm\Phi_-$ uniformly over $[-T/2,T/2]$ thus also in $L^2([-T/2,T/2])$.

\smallskip

Next, as $(\lambda_n)_{n\in \N}$ is increasing with $\lambda_0>0$ and
$\lambda_{n+1}-\lambda_n\to+\infty$, there exists $\alpha>0$ such that $\lambda_n\geq \alpha(n+1)$.
It follows that
$$
\sum_{k=0}^{+\infty}\frac{|a_k|}{\lambda_k}\leq \frac{1}{\alpha}\sum_{k=0}^{+\infty}\frac{|a_k|}{k+1}<+\infty
\quad\mbox{and}\quad
\sum_{k=0}^{+\infty}\frac{|a_k^{(n)}-a_k|}{\lambda_k}\leq \frac{1}{\alpha}\sum_{k=0}^{+\infty}\frac{|a_k^{(n)}-a_k|}{k+1}\to 0.
$$
As $|e^{2i\pi\lambda_kx}-1|\leq 2$, it follows that
$$
\Phi_+^{(n)}=\frac{1}{2i\pi}\sum_{k=K}^{K_n}\frac{a_k^{(n)}}{\lambda_k}\bigl(e^{2i\pi\lambda_kx}-1\bigr)
\to 
\Phi_+=\frac{1}{2i\pi}\sum_{k=K}^{\infty}\frac{a_k}{\lambda_k}\bigl(e^{2i\pi\lambda_kx}-1\bigr)
$$
where the series defining $\Phi_+$ is uniformly convergent over $[-T/2,T/2]$
and the convergence $\Phi_+^{(n)}\to\Phi_+$ is uniform over $[-T/2,T/2]$, thus also
in $L^2([-T/2,T/2])$. But we also know that $\Phi_+^{(n)}\to-\Phi_-$ in $L^2([-T/2,T/2])$
thus $\Phi_++\Phi_-=0$.

It remains to apply Kahane's extension of Ingham's Inequality to obtain that
\begin{eqnarray*}
0=\frac{1}{T}\int_{-T/2}^{T/2}|\Phi_+(t)+\Phi_-(t)|^2\,\mbox{d}t
&=&\frac{1}{T}\int_{-T/2}^{T/2}\abs{-\frac{1}{2i\pi}\sum_{k=0}^{+\infty}\frac{a_k}{\lambda_k}
+ \sum_{k=0}^{+\infty}\frac{a_k}{2i\pi\lambda_k}e^{2i\pi\lambda_k t}}^2\,\mbox{d}t\\
&\geq& A_2(T,\Lambda)\left(\left|\frac{1}{2i\pi}\sum_{k=0}^{+\infty}\frac{a_k}{\lambda_k}\right|^2
+\sum_{k=0}^{+\infty}\abs{\frac{a_k}{2i\pi\lambda_k}}^2\right)
\end{eqnarray*}
thus $a_k=0$ for all $k$ and we obtain the desired contradiction.

\medskip

The proof of \eqref{eq:inghamharl1} is similar, so we give less detail.
Elements of $\C^{2N+1}$ will be indexed from $-N$ to $N$, {\it i.e.}
$z=(z_{-N},\ldots,z_N)$ and will be considered as an element of $\C^{2M+1}$, $M\geq N$
and also as a sequence $(z_k)_{k\in\Z}$
by setting $z_k=0$ when $|k|>N$.

We again consider two norms on $\C^{2N+1}$, the $\ell^\infty$ norm
and (with a small abuse of notation)
$$
\normL{(a_{-N},\ldots,a_{N})}=\frac{1}{T}\int_{-T/2}^{T/2}\abs{\sum_{k=-N}^{N}a_ke^{2i\pi\lambda_k t}}\,\mbox{d}t.
$$
For every $N$ there exists $\tilde\kappa_N$ such that, for every $a\in\C^{2N+1}$,
$$
\tilde\kappa_N\norm{a}_\infty\leq\normL{a}.
$$
Ingham's Theorem asserts that one may choose $\kappa_N$ independently of $N$ provided $T$ is large enough.
Our aim is again to show that this is possible for every $T$ under our additional condition on $\lambda_k$.
Assume towards a contradiction that this is not possible. 

Then, for every integer $n\geq 1$, there exist
an integer $K_n\to+\infty$ with $K_{n+1}>K_n\geq K$
and $a^{(n)}=(a_{-K_n}^{(n)},\ldots,a_{K_n}^{(n)})\in\C^{2K_n+1}$ such that
$\|a^{(n)}\|_{\infty}=1$ while $\|a^{(n)}\|_{L^1_{K_n}}\leq\dfrac{1}{n}$. 
So, without loss of generality, we will assume that $K_{n+1}>K_n\geq K$ for every $n$.
Recall that we defined $K$ so that if $|k|\geq K$, $\lambda_{k+1}-\lambda_k\geq\dfrac{2}{T}$.

We split $a^{(n)}$ into two vectors
$$
a^{(n)}_-=(a_{-K+1}^{(n)},\ldots,a_{K-1}^{(n)})\in\C^{2K-1}
\quad\mbox{and}\quad
a^{(n)}_+=a^{(n)}-a^{(n)}_-.
$$
As $\|a^{(n)}_-\|_{\infty}\leq\|a^{(n)}\|_{\infty}= 1$,
there is no loss of generality in assuming that 
$$
a^{(n)}_-\to(a_{-K+1},\ldots,a_{K-1}).
$$

We again consider
$$
\ffi^{(n)}(t)=\sum_{k=-K_n}^{K_n}a_k^{(n)}e^{2i\pi\lambda_k t}\to 0
$$
in $L^1([-T/2,T/2])$,
$$
\ffi^{(n)}_-(t)=\sum_{k=-K+1}^{K-1}a_k^{(n)}e^{2i\pi\lambda_k t}
\to
\ffi_-(t)=\sum_{k=-K+1}^{K-1}a_ke^{2i\pi\lambda_k t}.
$$
in $L^1([-T/2,T/2])$ and
$$
\ffi^{(n)}_+=\ffi^{(n)}-\ffi^{(n)}_-=\sum_{K\leq|k|\leq K_n}a_k^{(n)}e^{2i\pi\lambda_k t}\to -\ffi_-
$$
in $L^1([-T/2,T/2])$. 

Using \eqref{eq:ingh} instead of \eqref{eq:HL2} we get, for $n\geq m$ 
$$
\frac{1}{T}\int_{-T/2}^{T/2}|\ffi^{(n)}_+(t)-\ffi^{(m)}_+(t)|\,\mbox{d}t
\geq A_1\left(T,\dfrac{2}{T}\right)\max_{K\leq|k|\leq K_n}\bigl|a_k^{(n)}-a_k^{(m)}\bigr|
$$
so that
$\bigl(a_k^{(n)}\bigr)_{k\in\Z}$ is a Cauchy sequence in $\ell^\infty$
and we call $a=(a_k)_{k\in\Z}$ its limit. Of course $\|a\|_\infty=1$
so that we will again reach a contradiction if we show that $a_k=0$ for all $k$.

To do so, we introduce again
$$
\Phi^{(n)}_\pm(x)=\int_0^x\ffi^{(n)}_\pm(t)\,\mbox{d}t
\quad\mbox{and}\quad
\Phi_-(x)=\int_0^x\ffi_-(t)\,\mbox{d}t
=\frac{1}{2i\pi}\sum_{k=-K+1}^{K-1}\frac{a_k}{\lambda_k}\bigl(e^{2i\pi\lambda_kx}-1\bigr)
$$
so that $\Phi^{(n)}_\pm\to\pm\Phi_-$ uniformly over $[-T/2,T/2]$ thus also in $L^2([-T/2,T/2])$.

\smallskip

Next, as $\lambda_k\not=0$ and $\sum_{k\in\Z}\dfrac{1}{1+|\lambda_n|}$ converges so is $\sum_{k\in\Z}\dfrac{1}{|\lambda_n|}$.
As $(a_k)\in\ell^\infty$ and $|a_k^{(n)}-a_k|\to 0$ in $\ell^\infty$
we get
$$
\sum_{k\in\Z}\frac{|a_k|}{|\lambda_k|}<+\infty
\quad\mbox{and}\quad
\sum_{k\in\Z}\frac{|a_k^{(n)}-a_k|}{|\lambda_k|}\to 0.
$$
As $|e^{2i\pi\lambda_kx}-1|\leq 2$, it follows that
$$
\Phi_+^{(n)}=\frac{1}{2i\pi}\sum_{K\leq|k|\leq K_n}\frac{a_k^{(n)}}{\lambda_k}\bigl(e^{2i\pi\lambda_kx}-1\bigr)
\to 
\Phi_+=\frac{1}{2i\pi}\sum_{|k|\geq K}\frac{a_k}{\lambda_k}\bigl(e^{2i\pi\lambda_kx}-1\bigr)
$$
where the series defining $\Phi_+$ is uniformly convergent over $[-T/2,T/2]$
and the convergence $\Phi_+^{(n)}\to\Phi_+$ is uniform over $[-T/2,T/2]$, thus also
in $L^2([-T/2,T/2])$. But we also know that $\Phi_+^{(n)}\to-\Phi_-$ in $L^2([-T/2,T/2])$
thus $\Phi_++\Phi_-=0$. Applying again Kahane's extension of Ingham's Inequality
we obtain 
$a_k=0$ for all $k$ which is the desired contradiction.

%
%
%
%
%

\section{1-periodic Schrödinger equation}

Recall that the Wiener algebra is defined as
$$
A(\T)=\{f\in L^1(\T)\,:\ \|f\|_{A(\T)}=\sum_{k\in \Z}|c_k(f)|<+\infty\}.
$$
\begin{theorem}
Let $u$ be a weak solution of the Schrödinger equation 
\begin{equation}
\begin{cases}
i\partial_t u(t,x)=\dfrac{1}{2\pi}\partial_x^2 u(t,x)&t\in\R,x\in\T\\
u_0=u(0,x)&x\in\R
\end{cases}\,,
\end{equation}
with initial value $u_0\in A(\T)$. Let $t_0\in\R$ and $x_0\in\T$.
Then
\begin{enumerate}
    \item For $a\in\R\setminus \Z$, for every $T>0$ there exists a constant $C(a,T)>0$ such that
    \begin{equation}\label{eq12}
        \frac{1}{T}\int_0^T |u(t_0+t,x_0+at)|\ud t \geqslant C\sum_{k\in \Z}\frac{|c_k(u_0)|}{1+|k|}.
    \end{equation}
    \item If $a \in \Z$, there exists $u_0\not=0$ such that $u(t_0+t,x_0+at)=0$
    for all $t$. In particular, \eqref{eq12} fails.
\end{enumerate}
\end{theorem}

\begin{remark}
Recall also that if $u_0\in H^s(\T)$ with $s>\dfrac{1}{2}$ then, with Cauchy-Schwarz,
$$
\|u_0\|_{A(\T)}\leq \left(\sum_{k\in\Z}\frac{1}{(1+|k|^2)^s}\right)^{1/2}
 \left(\sum_{k\in\Z}(1+|k|^2)^s|c_k(u_0)|^2\right)^{1/2}<+\infty.
$$
One may thus replace the condition $u_0\in A(\T)$ with a more familiar condition
like $u_0\in H^1(\T)$.
\end{remark}

\begin{proof} 
Write
$\dst
u_0(x)=\sum_{k\in \Z}c_ke^{2i\pi k x}
$
so that $\dst u(t,x)=\sum_{k\in \Z}c_ke^{2i\pi (k^2 t+k x)}$.
This series uniformly converges over $\R\times\T$ since $\sum|c_k|$ converges
thus $u$ is continuous. Further
$$
v(t)=u(t_0+t,x_0+at)=\sum_{k\in \Z}c_ke^{2i\pi k^2(t_0+t)+2i\pi k(x_0+at)}
=\sum_{k\in \Z}d_ke^{2i\pi \lambda_k t}
$$
with
$$
d_k=c_ke^{2i\pi(k^2t_0+kx_0)}\quad\mbox{and}\quad\lambda_k=k^2+ak.
$$
Note that $|d_k|=|c_k|$. On the other hand
\begin{eqnarray}
\lambda_k-\lambda_m&=&k^2+2ak-(m^2+am)=k^2-m^2+a(k-m)\notag\\
&=&(k-m)(k+m+a)\label{eq:lamdiff}.
\end{eqnarray}

Assume first that $a\in \Z$. This case was already treated in \cite{JK} but let us reproduce the proof here for completeness. In this case, the frequencies $(\lambda_k)$ satisfy the symmetry
property $\lambda_k=\lambda_{-a-k}$. Now fix $k\not=-a$ and notice that $-a-k\not=0$
so that, if we fix $c_k\not=0$
we can choose $c_{-a-k}$ so that $d_{-a-k}=-d_k$ that is
$$
c_{-a-k}=-c_ke^{2i\pi\bigl((k^2t_0+kx_0)-((-a-k)^2t_0+(-a-k)x_0)\bigr)}
=-c_ke^{-2i\pi\bigl(a(a+2k)t_0+x_0)\bigr)}.
$$
Setting
$$
u_0(x)=c_k\bigl(e^{2i\pi kt}-e^{-2i\pi\bigl(a(a+2k)t_0+2ax_0)\bigr)}e^{-2i\pi(a+k)t}\bigr)
$$
we obtain $u(t_0+t,x_0+at)=0$.

\smallskip

From now on, we assume that $a\notin\Z$ so that from \eqref{eq:lamdiff} we deduce that
$\lambda_k\not=\lambda_m$ when $k\not=m$. It will be convenient to write $a=2b$.
We can then further split the sequence $(\lambda_k)_{k\in\Z}$ into a disjoint
union, $(\lambda_k)_{k\in\Z}=(\lambda_k^+)_{k\geq 0}\cup
(\lambda_k^-)_{k\geq 1}$ with
$$
\lambda_k^+:=\lambda_{-[b]+k}= (-[b]+k)^2+2b(-[b]+k) \quad \text{for} \quad k\geqslant 0
$$
and
$$
\lambda_k^- :=\lambda_{-[b]-k}= (-[b]-k)^2+2b(-[b]-k) \quad \text{for}\quad k\geqslant 1.
$$
By definition 
$$
\lambda_0^+= [b]^2-2b[b] \quad \text{and} \quad \lambda_1^-=\lambda_0^++1-2(b-[b]).
$$
We will now distinguish two cases:

\smallskip

\noindent{\bf First case:} Assume that $\dfrac{1}{2}<b-[b]<1$ so that $\lambda_1^- <\lambda_0^+$.

\smallskip

In this case, the frequencies interlace as follows:
$$
\lambda_{k+1}^-<\lambda_k^+<\lambda_{k+2}^- \quad \text{for all}\quad k \geqslant 0.
$$
Indeed, for all $k\geqslant 0$
$$
\lambda_k^+-\lambda_{k+1}^-=2(2k+1)\left(b-[b]-\frac{1}{2}\right)>0 \quad \text{and} \quad \to \infty \quad \text{as} \quad k\to \infty
$$
and
$$
\lambda_{k+2}^--\lambda_k^+=4(k+1)\bigl(1-(b-[b])\bigr) >0 \quad \text{and} \quad \to \infty \quad \text{as} \quad k\to \infty
$$
with our hypothesis on $b-[b]$. In particular, if we set
$\mu_{2k}=\lambda_{k+1}^-=\lambda_{-[b]-k-1}$ and $\mu_{2k+1}=\lambda_{k}^+=\lambda_{[-b]+k}$ for $k\geq 0$ then
$0<\mu_{2k+1}-\mu_{2k}\to+\infty$.
Thus, from Theorem \ref{t11}, we get that
$$
\int_0^T|v(t)|\ud t\geqslant
 C(T)\sum_{k=0}^\infty\left(\frac{|c_{-[a]-k-1}|}{2k+1}+\frac{|c_{-[a]+k}|}{2k+2}\right).
$$
Finally, for $k\geq 0$, $2k+1\leq\alpha_a( |-[a]-k-1|+1)$ and $2k+2\leq\alpha_a(|-[a]+k|+1)$
with a constant $\alpha_a$ depending on $a$ only, so that
$$
\int_0^T|v(t)|\ud t\geqslant \frac{C(T)}{\alpha_a}\sum_{k\in\Z}\frac{|c_k|}{|k|+1}
$$
as claimed.

\smallskip

\noindent{\bf Second case:} $0<b-[b]<\dfrac{1}{2}$.

In this case, similar computations show that the frequencies interlace as
$$
\lambda_k^+<\lambda_{k+1}^-<\lambda_{k+1}^+ \quad \text{for all} \quad k\geqslant0
$$
with $\lambda_{k+1}^--\lambda_k^+,\lambda_{k+1}^+-\lambda_{k+1}^-\to+\infty$.
The remaining of the proof is the same and is thus omitted.

Note that $b-[b]\not=0,\dfrac{1}{2},1$ since $a=2b \notin \Z$ so all cases are now covered.
\end{proof}

\section{General Case}
Let $x$ in $\T=\R/\Z$ and $t \in \R^+$, we consider the following the equations
\begin{equation}
\label{eq:sgen}
\begin{cases}
i\partial_t u(t,x)=2\pi P\left(\dfrac{\partial_x}{2i\pi}\right)u\\
u_0=u(0,x)=\displaystyle\sum_{k\in \Z}c_ke^{2i\pi k x} \in \A(\T)
\end{cases}\,.
\end{equation}
where 
$$
P(X)=a_{n}X^{n}+a_{n-1}X^{n-1}+\ldots+a_1X+a_0
$$
with $n\geq 2$ and $a_n\not=0$. There is no
loss of generality in assuming that $a_n>0$.

If $u_0(x)=\dst\sum_{k\in\Z}c_ke^{2i\pi kx}\in A(\T)$, then the solution to this system is given by
$$
u(t,x)=\sum_{k\in \Z}c_ke^{-2i\pi P(k)t}e^{2i\pi kx}.
$$
Again, this is a continuous function.

Let $a\in \R$ to be chosen later. For any $(t_0,x_0)\in \R\times \T$, we define 
$$
u(t_0+t,x_0+at)=\sum_{k\in \Z}c_ke^{-2i\pi P(k)(t_0+t)+2i\pi k(x_0+at)}\\
=\sum_{k\in \Z}d_ke^{-2i\pi \lambda_k t}
$$
with
$$
d_k=c_ke^{-2i\pi (P(k)t_0-kx_0)}
\quad\mbox{and}\quad
\lambda_k=P(k)-ak.
$$
Note that $\lambda_k-\lambda_m=(k-m)\bigl(Q(k,m)-a\bigr)$
with
\begin{eqnarray*}
Q(k,m)&=&a_{n}(k^{n-1}+k^{n-2}m+\ldots+m^{n-1})+a_{n-1}(k^{n-2}+\ldots+m^{n-2})+\ldots+a_1\\
&=&\sum_{\ell=1}^n a_\ell\sum_{j=0}^{\ell-1} k^{\ell-1-j}m^j.
\end{eqnarray*}

Define 
$$
E=\{Q(k,m),~k,m \in \Z~\text{such that}~k\neq m\}
$$
which is countable (thus of measure $0$). 

\begin{theorem}
Let $u$ be any solution of the Schrödinger equation \eqref{eq:sgen} with initial value $u_0=\dst\sum_{k\in\Z}c_ke^{2i\pi kx}\in A(\T)$. Then
\begin{enumerate}
    \item If $a\notin E$, for all $T>0$ there exists a constant $C(a,T)>0$ such that
\begin{equation}\label{e22even}
        \frac{1}{T}\int_0^T |u(t_0+t,x_0+at)|\ud t \geqslant C\max_{k\in \Z}|c_k|.
    \end{equation}
    If $n$ is even, there also exists a constant $C(a,T)>0$ such that
    \begin{equation}\label{e22odd}
        \frac{1}{T}\int_0^T |u(t_0+t,x_0+at)|\ud t \geqslant C\sum_{k\in \Z}\frac{|c_k|}{1+|k|}.
    \end{equation}
    \item If $a\in E$ then both \eqref{e22even}-\eqref{e22odd} fail.
\end{enumerate}
\end{theorem}

An $L^2$ analogue of this result can be found in \cite{WW}. 

\begin{proof}   
The last part of the theorem is the same as for the Schrödinger equation in the previous section.
Indeed, if $a\in E$, we can choose two indexes $k\not=m$ such that $\lambda_k=\lambda_m$ and then choose $c_k,c_m$
such that $d_k=-d_m$. Taking $u_0=c_ke^{2i\pi kt}+c_me^{2i\pi mt}$,
the corresponding solution $u$ satisfies $u(t_0+t,x_0+at)=0$.

\medskip

We now assume that $a\notin E$ so that $\lambda_k\not=\lambda_m$ for all $k,m\in\Z$.
We will further show that the $(\lambda_k)'s$ can be ordered as a sequence with gaps going to infinity.
Here we need to distinguish between $n$ even or odd. We start with the odd case.

\medskip

If $n$ is odd, then $\lambda_k=P(k)-ak\to\pm\infty$ when $k\to\pm\infty$. Note also that, as $P$ has degree at least $3$, $\sum\dfrac{1}{1+|\lambda_k|}$ converges.

Further $\lambda_{k+1}-\lambda_k=Q(k+1,k)-a=a_nk^{n-1}+o(k^{n-1})\to+\infty$ when $k\to\pm\infty$.
so that, there exists $K$ such that, for $k\geq K$, $\lambda_k$ is increasing as well as for $k\leq -K$.
There further exists $K'\geq K$ such that, if $k,\ell\geq K'$,
then 
$$
\lambda_{-\ell}\leq \min_{|j|\leq K}\lambda_j\leq\max_{|j|\leq K}\lambda_j\leq\lambda_k.
$$
We then define $(\mu_k)_{|k|\leq K'}$ as an ordering of $(\lambda_k)_{|k|\leq K'}$
and $\mu_k=\lambda_k$ for $|k|>K'$. Note that those $\lambda_k$'s are not one of the $(\mu_k)_{|k|\leq K'}$'s.
Then $(\mu_k)_{k\in\Z}$ is an increasing sequence
with gaps $\mu_{k+1}-\mu_k\to+\infty$ when $k\to\pm\infty$. We can then apply \eqref{eq:inghamharl1}
to conclude.

\medskip

We now assume that $n=2p$ is even. In this case $\lambda_k=P(k)-ak\to+\infty$ when $k\to\pm\infty$
and $\lambda_{k+1}-\lambda_k\to\pm\infty$ when $k\to\pm\infty$.
In this case, the ordering needs to be made differently.

The idea is rather simple, there is an oscillating part and we are going to show that, for $k,\ell$ large,
the $\lambda_k$'s and $\lambda_{-\ell}$'s interlace. In the generic case
we actually have $\lambda_{k+q_0}<\lambda_{-k}<\lambda_{k+q_0+1}$ for some fixed $q_0$
and large enough $k$. This shows that, for some $K_0$, $(\lambda_k)_{k\notin\{-K_0,\ldots,K_0+q_0\}}$
can be rearranged in an increasing way as $\lambda_{K_0+q_0+1},\lambda_{-K_0-1},\lambda_{K_0+q_0+2},\lambda_{-K_0-2},\ldots$. The finite number of remaining $\lambda_k$'s are rearranged separately and, provided $K_0$
is large enough, they can be put at the start and the resulting sequence $(\mu_k)_{k\geq 0}$ is then increasing with gaps going to infinity. A key aspect of this construction is that each $\mu_k$ is a $\lambda_{k'}$
with $\bigl|k-|k'|\bigr|\leq C_\Lambda$ depending only on $\Lambda$.
The idea is the same in the exceptional case.

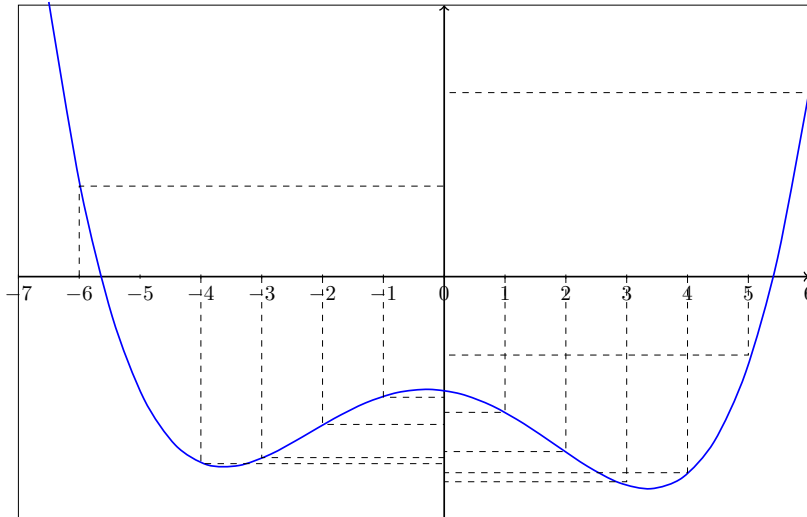
\begin{center}
\scalebox{0.8}{
\begin{tikzpicture}
\draw (-7,-4) -- (-7,4.5) -- (6,4.5) -- (6,-4) -- (-7,-4);
\draw[thick,->] (-7,0) -- (6,0);
\draw[thick,->] (0,-4) -- (0,4.5);
\foreach \x in {-7,-6,-5,-4,-3,-2,-1,0,1,2,3,4,5,6}
    \draw(\x,1pt)--(\x,-1pt) node[below]{$\x$};
\draw[blue, thick,smooth] [domain=-6.5:6] plot(\x,{(\x/1.5-3.2)*(\x/1.5-0.3)*(\x/1.5+0.7)*(\x/1.5+3.3)/20-2});
\draw[very thin,dashed] (-1,0) -- (-1,-2) -- (0,-2);
\draw[very thin,dashed] (-2,0) -- (-2,-2.45) -- (0,-2.45);
\draw[very thin,dashed] (-3,0) -- (-3,-3) -- (0,-3);
\draw[very thin,dashed] (-4,0) -- (-4,-3.1) -- (0,-3.1);
\draw[very thin,dashed] (-6,0) -- (-6,1.5) -- (0,1.5);
\draw[very thin,dashed] (6,0) -- (6,3.05) -- (0,3.05);
\draw[very thin,dashed] (5,0) -- (5,-1.3) -- (0,-1.3);
\draw[very thin,dashed] (4,0) -- (4,-3.25) -- (0,-3.25);
\draw[very thin,dashed] (3,0) -- (3,-3.4) -- (0,-3.4);
\draw[very thin,dashed] (2,0) -- (2,-2.9) -- (0,-2.9);
\draw[very thin,dashed] (1,0) -- (1,-2.25) -- (0,-2.25);
\end{tikzpicture}
}
\begin{minipage}{0.9\textwidth}
\small
The picture shows the case of a polynomial $P$ of degree $4$.
The reordering here is $\mu_0=\lambda_3$, $\mu_1=\lambda_4$,
$\mu_2=\lambda_{-4}$, $\mu_3=\lambda_{-3}$
$\mu_4=\lambda_2$, $\mu_5=\lambda{-2}$, $\mu_6=\lambda_1$,
$\mu_7=\lambda_1$, $\mu_8=\lambda_{-5}$ (not represented to keep the picture readable), $\mu_9=\lambda_{-1}$,
$\mu_{10}=\lambda_0$, $\mu_{11}=\lambda_5$, $\mu_{12}=\lambda_{-6}$,
$\mu_{13}=\lambda_6$ and more generally $\mu_{12+2k}=\lambda_{-6-k}$
while $\mu_{13+2k}=\lambda_{6+k}$.
\end{minipage}
\end{center}

Let us now be more precise.

We again take $K$
such that from $-\infty$ to $K$,
$\lambda_k$ is decreasing while from $K$ to $+\infty$, $\lambda_k$ is increasing
and define $K'$ such that if $k,\ell\geq K'$,
then 
$$
\max_{|j|\leq K}\lambda_j\leq\lambda_{-\ell},\lambda_k.
$$

Next, an easy computation shows that
$$
\sum_{j=0}^\ell(-1)^j=\begin{cases}1&\mbox{if $\ell$ is even}\\
0&\mbox{if $\ell$ is odd}\end{cases}
\quad\mbox{and}\quad
\sum_{j=0}^\ell(-1)^jj=\begin{cases}\ell/2&\mbox{if $\ell$ is even}\\
-(\ell+1)/2&\mbox{if $\ell$ is odd}\end{cases}
$$
so that
$$
\lambda_{k+q}-\lambda_{-k}=
(2k+q)\left(\sum_{\ell=1}^{2p}a_\ell\sum_{j=0}^{\ell-1}(-1)^j(k+q)^{\ell-1-j}k^j-a\right).
$$
But
\begin{eqnarray*}
a_{2p}\sum_{j=0}^{2p-1}(-1)^j(k+q)^{\ell-1-j}k^j&=&a_{2p}k^{2p-1}\sum_{j=0}^{2p-1}(-1)^j
+a_{2p}qk^{2p-2}\sum_{j=0}^{2p-2}(-1)^jj+o(k^{2p-2})\\
&=&(p-1)a_{2p}qk^{2p-2}+o(k^{2p-2})
\end{eqnarray*}
and
\begin{eqnarray*}
a_{2p-1}\sum_{j=0}^{2p-2}(-1)^j(k+q)^{2p-2-j}k^j
&=&a_{2p-1}k^{2p-2}+o(k^{2p-2})
\end{eqnarray*}
so that
\begin{equation}
\label{eq:leading}
\lambda_{k+q}-\lambda_{-k}=\begin{cases}\bigl((p-1)a_{2p}q+a_{2p-1}\bigr)k^{2p-1}+o(k^{2p-1})&\mbox{if }p\geq 2\\
\bigl(a_{2}q+a_1-a\bigr)k+o(k)&\mbox{if }p=1
\end{cases}.
\end{equation}

Set $\alpha_q=\begin{cases}a_{2}q+a_1-a&\mbox{if }p=1\\ (p-1)a_{2p}q+a_{2p-1}&\mbox{if }p\geq2\end{cases}$
so that $\lambda_{k+q}-\lambda_{-k}=\alpha_qk^{2p-1}+o(k^{2p-1})$.

There are now two cases:

\medskip

\noindent{\bf Case 2.1} {\em For every $q$, $\alpha_q\not=0$}

\medskip

Then there exists
$q_0$ such that $\alpha_{q_0}>0$ and $\alpha_{q_0-1}<0$. But then, $\lambda_{k+q_0}-\lambda_{-k}\to+\infty$
while $\lambda_{k+q_0-1}-\lambda_{-k}\to-\infty$.

We now take $K''>\max(K'-q_0,K')$ such that, for $k\geq K''$, $\lambda_{k+q_0}-\lambda_{-k}>0$
and $\lambda_{k+q_0-1}-\lambda_{-k}<0$, that is $\lambda_{k+q_0-1}<\lambda_{-k}<\lambda_{k+q_0}$.
The choice of $K''$ also implies that $\lambda_{-K''+1},\ldots, \lambda_{K''+q_0-1}$
are all $<\min(\lambda_{-K''},\lambda_{K''+q_0})$.
We can thus reorder those terms as an increasing sequence $(\mu_k)_{k=0,\cdots,\hat K}$
with $\hat K=2K''+q_0-2$,
that we then complete into a sequence $(\mu_k)_{k\in\N}$ by adding successively
a term $\lambda_{K''+k+q_0}$ and a term $\lambda_{-K''-k}$ and the resulting sequence is an increasing
rearrangement of $(\lambda_k)$ such that $\mu_k\to+\infty$ and $\mu_{k+1}-\mu_k\to+\infty$.
Note that if we define $\sigma$ the mapping $\Z\to\N$ defined by $\mu_k=\lambda_{\sigma(k)}$
then there is a constant $C_\Lambda$ such that $\bigl||k|-\sigma(k)\bigr|\leq C_\Lambda$.

It follows from \eqref{eq:harrauxl1} that
\begin{eqnarray*}
\frac{1}{T}\int_{-T/2}^{T/2}|u(t_0+t,x_0+at)|\,\mbox{d}t&=&
\frac{1}{T}\int_{-T/2}^{T/2}\abs{\sum_{k\in \Z}d_ke^{-2i\pi \lambda_k t}}\,\mbox{d}t\\
&=&\frac{1}{T}\int_{-T/2}^{T/2}\abs{\sum_{k=0}^{+\infty}d_{\sigma^{-1}(k)}e^{-2i\pi \mu_k t}}\,\mbox{d}t\\
&\geq&\tilde A_1(T,\Lambda)\sum_{k=0}^{+\infty}\frac{|d_{\sigma^{-1}(k)}|}{1+k}
=\tilde A_1(T,\Lambda)\sum_{j\in\Z}\frac{|c_j|}{1+\sigma(j)}\\
&\geq&\frac{\tilde A_1(T,\Lambda)}{1+C_\Lambda}\sum_{k\in\Z}\frac{|c_k|}{1+|k|}.
\end{eqnarray*}
Note that the series $\sum_{k\in \Z}d_ke^{-2i\pi \lambda_k t}$ is uniformly
convergent so that we can re-order terms.

\medskip

\noindent{\bf Case 2.2} {\em There exists $q_0$, such that $\alpha_{q_0}=0$.}

\medskip

The proof is essentially the same, but the interlacing of the $\lambda_k$ and $\lambda_{-\ell}$
for $k,\ell$ large may be different. This comes from the fact that the leading term in \eqref{eq:leading} is now $0$.
Nevertheless, $\alpha_{q_0+1}>0$ and $\alpha_{q_0-1}<0$ so that, for $k$ large enough 
$\lambda_{k+q_0+1}-\lambda_{-k}>0$ while $\lambda_{k+q_0-1}-\lambda_{-k}<0$. 
So, for each $k$, either $\lambda_{k+q_0-1}<\lambda_{k+q_0}<\lambda_{-k}<\lambda_{k+q_0+1}$
or $\lambda_{k+q_0-1}<\lambda_{-k}<\lambda_{k+q_0}<\lambda_{k+q_0+1}$ (actually
only one can occur for $k$ large enough) and we define the rearrangement $\mu_k$
accordingly.
\end{proof}

\end{document}